\documentclass[12pt]{amsart}

\usepackage{amscd}
\usepackage{srcltx}
\usepackage{amsmath}
\usepackage{amsfonts}
\usepackage{verbatim}
\usepackage{amsthm}
\usepackage{amssymb, mathrsfs}
\usepackage{amscd}
\usepackage{latexsym}
\usepackage{array}
\usepackage{enumerate}
\usepackage{longtable}
\usepackage[all]{xypic}
\usepackage[latin1]{inputenc}
\usepackage{epsfig}

\ifx\pdfpagewidth\undefined 
\usepackage[T1]{fontenc}
\else 
\usepackage[OT1]{fontenc} 
\fi 
\usepackage{amsfonts,amssymb,latexsym,longtable,amsmath}

\newtheorem{thm}{Theorem}[section]
\newcommand{\bthm}{\begin{thm}} \newcommand{\ethm}{\end{thm}}
\newtheorem{prop}[thm]{Proposition}
\newcommand{\bprp}{\begin{prop}} \newcommand{\eprp}{\end{prop}}
\newtheorem{fact}[thm]{Fact}
\newcommand{\bfct}{\begin{fact}} \newcommand{\efct}{\end{fact}}
\newtheorem{prob}[thm]{Problem}
\newcommand{\bprb}{\begin{prob}} \newcommand{\eprb}{\end{prob}}
\newtheorem{lem}[thm]{Lemma}
\newcommand{\blem}{\begin{lem}} \newcommand{\elem}{\end{lem}}
\newtheorem{claim}[thm]{Claim}
\newcommand{\bclm}{\begin{claim}} \newcommand{\eclm}{\end{claim}}
\newtheorem{cor}[thm]{Corollary}
\newcommand{\bcor}{\begin{cor}} \newcommand{\ecor}{\end{cor}}
\newtheorem{conj}[thm]{Conjecture}
\newcommand{\bcnj}{\begin{conj}} \newcommand{\ecnj}{\end{conj}}
\theoremstyle{definition}
\newtheorem{defn}[thm]{Definition}
\newcommand{\bdfn}{\begin{defn}} \newcommand{\edfn}{\end{defn}}
\newtheorem{spec}[thm]{Specializing}
\newcommand{\bspc}{\begin{spec}} \newcommand{\espc}{\end{spec}}
\theoremstyle{remark}
\newtheorem{rem}[thm]{Remark}
\newcommand{\brem}{\begin{rem}} \newcommand{\erem}{\end{rem}}
\newtheorem{cnv}[thm]{Convention}
\newcommand{\bcnv}{\begin{cnv}} \newcommand{\ecnv}{\end{cnv}}
\newtheorem{exam}[thm]{Example}
\newcommand{\bexm}{\begin{exam}} \newcommand{\eexm}{\end{exam}}
\newcommand{\bpf}{\begin{proof}} \newcommand{\epf}{\end{proof}}

\newcommand{\R}{\mathbb R}

\newcommand{\Z}{\mathbb Z}
\newcommand{\T}{\mathbb T}
\newcommand{\N}{\mathbb N}

\renewcommand{\phi}{\varphi}
\renewcommand{\theta}{\vartheta}

\begin{document}
\title{On the structure of abelian profinite groups}

\author[M. Ferrer]{M. Ferrer}
\address{Universitat Jaume I, Instituto de Matem\'aticas de Castell\'on,
Campus de Riu Sec, 12071 Castell\'{o}n, Spain.}
\email{mferrer@mat.uji.es}
\author[S. Hern\'andez]{S. Hern\'andez}
\address{Universitat Jaume I, Departamento de Matem\'{a}ticas,
Campus de Riu Sec, 12071 Castell\'{o}n, Spain.}
\email{hernande@mat.uji.es}

\thanks{ The first-listed
author acknowledges partial support by the Universitat Jaume I, grant P1171B2015-77. 
The second-listed author acknowledges partial financial support by the Spanish Ministry of
Science, grant MTM2008-04599/MTM; and the Universitat Jaume I, grant P1171B2015-77}

\date{November/19/2018}

\maketitle \setlength{\baselineskip}{24pt}
\setlength{\parindent}{1cm}
\begin{abstract}
A subgroup $G$ of a product $\prod\limits_{i\in\mathbb{N}}G_i$ is \emph{rectangular} if there are subgroups $H_i$ of $G_i$
such that $G=\prod\limits_{i\in\mathbb{N}}H_i$. We say that $G$ is \emph{weakly rectangular} if there are finite subsets
$F_i\subseteq \mathbb{N}$ and subgroups $H_i$ of $\bigoplus\limits_{j\in F_i} G_j$ that satisfy $G=\prod\limits_{i\in\mathbb{N}}H_i$.
In this paper we discuss when a closed subgroup of a product is
weakly rectangular. Some possible applications to the theory of group codes are also highlighted.
\end{abstract}
\vspace{1cm}

\section{Introduction}

For a family $\{G_i : i\in \mathbb{N} \}$ of topological groups, let $\bigoplus\limits_{i\in \mathbb{N}} G_i$
denote the subgroup of elements $(g_i)$ in the product $\prod\limits_{i\in\mathbb{N}}G_i$
such that $g_i=e$ for all but finitely many indices $i\in \mathbb{N}$. 
A subgroup  $G\leq \prod\limits_{i\in\mathbb{N}}G_i$ is called \emph{weakly controllable} if
$G\cap \bigoplus\limits_{i\in \mathbb{N}} G_i$ is dense in $G$, that is, if $G$ is generated by
its elements of finite support. The group $G$ is called \emph{weakly observable} if
$G\cap \bigoplus\limits_{i\in \mathbb{N}} G_i=\overline G\cap \bigoplus\limits_{i\in \mathbb{N}} G_i$,
where $\overline G$ stands for the closure of $G$ in $\prod\limits_{i\in\mathbb{N}}G_i$ for
the product topology. Although the notion of (weak) controllability was coined by Fagnani earlier
in a broader context (cf. \cite{fagnani:adv97}), both notions
were introduced in the area of of coding theory by Forney and Trott (cf. \cite{forney_trott:04}).
They observed that if
the groups $G_i$ are locally compact abelian, then controllability and observability are dual properties
with respect to the Pontryagin duality: If $G$ is a closed subgroup of $\prod\limits_{i\in\mathbb{N}}G_i$,
then it is weakly controllable if and only if its annihilator
$G^\bot=\{\chi\in \widehat{\prod\limits_{i\in\mathbb{N}}G_i} : \chi(G)=\{0\}\}$ is a weakly observable subgroup
of $\prod\limits_{i\in\mathbb{N}}\widehat{G_i}$ (cf. \cite[4.8]{forney_trott:04}).

In connection with the properties described above, the following definitions was introduced in \cite{GL:2012}.

\bdfn
A subgroup $G$ of a product $\prod\limits_{i\in\mathbb{N}}G_i$ is \emph{rectangular} if there are subgroups $H_i$ of $G_i$
such that $G=\prod\limits_{i\in\mathbb{N}}H_i$. We say that $G$ is \emph{weakly rectangular} if there are finite subsets
$F_i\subseteq \mathbb{N}$ and subgroups $H_i$ of $\bigoplus\limits_{j\in F_i} G_j$ that satisfy $G=\prod\limits_{i\in\mathbb{N}}H_i$.
We say that $G$ is a \emph{subdirect product} of the family $\{G_i\}_{i\in I}$ if $G$ is weakly rectangular and
$G\cap\bigoplus\limits_{i\in I} G_i=\bigoplus\limits_{i\in\mathbb{N}}H_i$.
\edfn
\medskip

\noindent The observations below are easily verified. 

\begin{enumerate}
\item{Weakly rectangular subgroups and rectangular subgroups of $\prod\limits_{i\in\mathbb{N}}G_i$ are weakly controllable.}
\item{If each $G_i$ is a pro-$p_i$-group for some prime $p_i$, and all $p_i$ are distinct, then every closed
subgroup of the product $\prod\limits_{i\in\mathbb{N}}G_i$ is rectangular, and thus is a subdirect product.}
\item{If each $G_i$ is a finite simple non-abelian group and all $G_i$ are distinct, then every closed normal subgroup of the
product  $\prod\limits_{i\in\mathbb{N}}G_i$ is rectangular, and thus a subdirect product.}
\end{enumerate}
\medskip

The main goal addressed in this paper is to study to what extent the converse of these observations hold.
In particular, we are interested in the following question (cf. \cite{GL:2012}):

\bprb\label{Question}
Let $\{G_i : i\in \mathbb{N} \}$ be a family of compact metrizable groups, and $G$
a closed subgroup of the product $\prod\limits_{i\in\mathbb{N}}G_i$. If $G$ is weakly controllable, that is,
$G\cap \bigoplus\limits_{i\in \mathbb{N}} G_i$ is dense in $G$, what
can be said about the structure of $G$? In particular, under what additional conditions is the group
$G$ a subdirect product of $\{G_i : i\in \mathbb{N} \}$, that is, weakly rectangular and
$G\cap\bigoplus\limits_{i\in I} G_i=\bigoplus\limits_{i\in\mathbb{N}}H_i$, where each $H_i$ is a subgroup
of $\bigoplus\limits_{j\in F_i} G_j$ for some $F_i\subseteq \mathbb{N}$?
\eprb
\medskip

In connection with this question, the following result was established in \cite{FHS:2017}.

\bthm\label{DHS:products}
Let $I$ be a countable set, $\{G_i:i\in I\}$ be a family of finite abelian groups
and $G=\prod_{i\in I} G_i$ be its direct product.
Then every closed weakly controllable subgroup $H$ of $G$
is topologically isomorphic to a direct product
of finite cyclic groups.
\ethm
\medskip

Unfortunately, this result does not answer Problem \ref{Question}, which remains open to the best of our knowledge.
Finally, it is pertinent to mention here that the relevance of these notions stem from coding theory where they appear
in connection with the study of convolutional group codes \cite{forney_trott:04,rosenthal}.
However, similar notions had been studied in symbolic dynamics previously. Thus, the notions
of weak controllability and weak observability are related to the concepts of \emph{irreducible shift} and
\emph{shift of finite type}, respectively, that appear in symbolic dynamics.
Here, we are concerned with abelian profinite groups and
our main interest is to clarify the overall topological and algebraic structure of abelian profinite groups
that satisfy any of the properties introduced above.
In the last section, we shall also consider some possible applications to the study of group codes.

\section{Basic facts}
\subsection{Pontryagin-van Kampen duality}

One of the main tools in this research is Kaplan's extension of Pontryagin van-Kampen duality
to infinite products of locally compact abelian (LCA) groups. In like manner that Pontryagin-van Kampen
duality has proven to be essential in understanding the structure of LCA, the extension accomplished
by Kaplan for cartesian products and direct sums \cite{kaplan1,kaplan2} (and some other subsequent results) have
established duality methods as a powerful tool outside the class of LCA groups and they have been
widely used in the study of group codes.


We recall the basic properties of topological
abelian groups and the celebrated Pontryagin-van Kampen duality.

Let $G$ be a commutative locally-compact group. A character $\chi$ of $G$ is a continuous homomorphism
$\chi : G\longrightarrow \T$
where $\T$ is the multiplicative group of complex numbers of modulus $1$.
The characters form a group $\widehat G$, called \emph{dual group},
which is given the topology of uniform convergence on compact subsets of $G$.
It turns out that $\widehat G$ is locally compact and there is a canonical 
evaluation homomorphism
$$\mathcal E_G : G \longrightarrow \widehat{\widehat G
}.$$
\bthm\label{th_Pontryagin-van Kampen}
The evaluation homomorphism $\mathcal E_G$ is an isomorphism of
topological groups.
\ethm

\noindent Some examples: $$\widehat \T\cong \Z,\ \widehat \Z \cong \T,\ \widehat \R \cong \R,\ \widehat{(\Z/n)}\cong \Z/n.$$
(some groups are self dual, such as finite abelian groups or
the additive group of the real numbers)

Pon\-trya\-gin-van Kampen dua\-li\-ty
establishes a dua\-li\-ty between the subcategories of compact
and discrete abelian groups.
If $G$ denotes a compact abelian group and $\Gamma $
denotes its dual group, we have the following equivalences between
topological properties of $G$ and algebraic properties of $\Gamma$:

\begin{enumerate}{}
\item[(i)]$\quad weight(G)=|\Gamma |\quad $(metrizablity$\Leftrightarrow
|\Gamma |\leq \omega$);

\item[(ii)]\quad $G$ is connected $\Leftrightarrow $ $\Gamma $ is torsion
free;

\item[(iii)]\quad $Dim(G)=0\Leftrightarrow $ $\Gamma $ is torsion; and

\item[(iv)]\quad $G$ is monothetic $\Leftrightarrow $ $\Gamma $ is
isomorphic to a subgroup of $\Bbb{T}_{d}$.
\end{enumerate}

In general, it is said that a topological abelian group $(G,\tau )$ satisfies
\emph{Pontryagin duality} (is \emph{P-reflexive} for short,) if the evaluation map
$\mathcal E_{G}$ is a topological isomorphism onto. We refer to the survey by
Dikranjan and Stoyanov \cite{ds:duality} and the  monographs by
Dikranjan, Prodanov and Stoyanov \cite{dps:duality} and Hofmann, Morris \cite{hm:compact}
in order to find the basic results about Pontryagin-van Kampen duality.

The following result, due to Kaplan \cite{kaplan1,kaplan2} is essential in the applications
of duality methods.
\bthm[Kaplan]\label{th Kaplan}
Let $\{G_{i}:i\in I\}$ be a family of $P$-reflexive groups. Then:

\begin{enumerate}
\item[(i)]\qquad The direct product $\prod\limits_{i\in I}G_{i}$ is $P$-reflexive.

\item[(ii)]\qquad The direct sum $\bigoplus\limits_{i\in I}G_{i}$ equipped with a
suitable topology is a $P$-reflexive group.

\item[(iii)]\qquad It holds:
\begin{eqnarray*}
(\prod\limits_{i\in I}G_{i})^{\widehat{}} &\cong &\bigoplus_{i\in I}\widehat{%
G_{i}} \\
(\bigoplus\limits_{i\in I}G_{i})^{\widehat{}} &\cong &\prod_{i\in I}\widehat{%
G_{i}}
\end{eqnarray*}
\end{enumerate}
\ethm

Kaplan also set the problem of characterizing the class of topological Abelian
groups for which Pontryagin duality holds.

Let $g\in G$ and $\chi\in\widehat{G}$, it is said that $g$ and $\chi$
are \emph{orthogonal} when $\langle g,\chi\rangle=1$. Given
$S\subseteq G$ and $S_1\subseteq \widehat G$ we define \emph{the
orthogonal (or annihilator)}  of  $S$ and $S_1$ as
$$
S^{\perp}=\{\chi\in\widehat{G}\,:\,\langle g,\chi\rangle=1\;
\forall g\in S\}$$
and
$$S_1^{\perp}=\{g\in
G\,:\,\langle g,\chi\rangle=1\; \forall \chi\in S_1\}.
$$
\noindent Obviously $G^{\perp}=\{e_{\widehat{G}}\}$ and $\widehat{G}^{\perp}=\{e_{G}\}$.

The following result has also many applications in connection with duality theory.

\bthm
Let $S$ and $R$ be subgroups of a LCA group $G$ such that $S\leq
R\leq G$. Then we have $\widehat{R/S} \cong S^{\perp}/R^{\perp}$.
\ethm
\bcor
Let $H$ be a closed subgroup of a LCA group $G$.  Then
$\widehat{G/H} \cong H^{\perp}$ and $\widehat{H} \cong
\widehat{G}/H^{\perp}$.
\ecor


\subsection{Abelian profinite groups}

Our main results concern the structure of abelian profinite groups that appear in coding theory.
Firstly, we recall some basic definitions and terminology.
For every group $G$ let us denote by $(G)_p$ the largest $p$-subgroup of $G$ and
$\mathbb{P}_G=\{p\in\mathbb{P}\,:\, G\text{ contains a }p-\text{subgroup}\}$
where $p\in \mathbb{P}_G$ and $\mathbb{P}$ is the set of all prime numbers.
An element $g$ of a $p$-primary group $G$ is said to have \emph{ finite height in $G$} $h$ if this is
the largest natural number $n$ such that the equation $p^nx = g$ has a solution $x\in G$.
We say that $g$ has \emph{infinite height} if  the solution exists  for all $n\in \N$.
Here on, the symbol $G[p]$ denotes the subgroup consisting of all
elements of order $p$. It is well known that $G[p]$ is a vector space on the field $\Z(p)$.

\section{Order controllable groups}

\bdfn Let $\{G_i : i\in \mathbb{N} \}$ be a family of topological groups and $\mathcal{C}$ a subgroup of $\prod\limits_{i\in \mathbb{N}} G_i$. We have the
following notions:

$\mathcal{C}$ is \emph{weakly controllable} if $\mathcal{C}\bigcap\bigoplus\limits_{i\in
\mathbb{N}}G_i$ is dense in $\mathcal{C}$.

$\mathcal{C}$ is \emph{(uniformly) controllable} if for every $i\in \mathbb{N}$ there is $n_i\in \mathbb{N}$ such that if $c\in\mathcal{C}$ there
exists $c_1\in\mathcal{C}$ such that $c_1|_{[1,i]}=c|_{[1,i]}$ and $c_1|_{]n_i,+\infty[}=0$
(we assume that $n_i$ is the less natural number satisfying this property).
This implies that there exists $c_2\in \mathcal{C}$ such that $c=c_1+c_2$, $supp (c_1)\subseteq [1,n_i]$ and $supp(c_2)\subseteq [i+1,+\infty[$.

$\mathcal{C}$ is \emph{order-controllable} if for every $i\in \mathbb{N}$ there is $n_i\in \mathbb{N}$ such that if $c\in\mathcal{C}$ there
exist $c_1$ and $c_2$ in $\mathcal{C}$ such that $c=c_1+c_2$, with
$supp (c_1)\subseteq [1,n_i]$, $supp (c_2)\subseteq [i+1,+\infty[$ and $order(c_1)\leq order(c|_{[1,n_i]})$
(again, we assume that $n_i$ is the less natural number satisfying this property).
As a consequence, we also have that $order(c_2)\leq order(c)$. Here, the order $c$ is taken in the usual sense, considering $c$ as an element of the group $\mathcal C$.

\edfn

Every controllable group is weakly controllable and, if the groups $G_i$ are finite, then the notions of controllability and weakly controllability are equivalent
(see \cite{FHS:2017}). The following result partially answers Problem \ref{Question} for $p$-groups. The proof can be founded in \cite{FH:2018}.

\bthm\label{th_products_p}
Let $\{G_i : i\in \mathbb{N} \}$ be a family of  finite, abelian, $p$-groups and let $G=\prod\limits_{i\in \mathbb{N}} G_i$.
If $\mathcal{C}$ is an infinite closed subgroup of $G$ which is order-controllable, then
$\mathcal{C}$ is weakly rectangular. In particular, there is a sequence $\{y_m: m\in \mathbb{N} \}\subseteq \mathcal{C}\bigcap\bigoplus\limits_{i\in \mathbb{N}}G_i$ such that $\mathcal{C}$ is
topologically isomorphic to $\prod\limits_{m\in \mathbb{N}} <y_m>$.
\ethm

This result extends directly to general products of finite abelian groups and gives a partial answer to Question \ref{Question}.

\bthm\label{th_products}
Let $\mathcal{C}$ be an order-controllable, closed, subgroup of a countable product $G=\prod\limits_{i\in\mathbb{N}}G_i$ of finite abelian groups $G_i$.
Then $\mathcal{C}$ is weakly rectangular. In particular,
there is a sequence $\{y_{m}^{(p)}\,:\,m\in\mathbb{N},p\in\mathbb{P}_G\}\subseteq G\cap(\bigoplus\limits_{i\in\mathbb{N}}G_i)$ such that
$\{y_m^{(p)}\,:\,m\in\mathbb{N}\}\subseteq(G\cap(\bigoplus\limits_{i\in\mathbb{N}}G_i))_p$
and $G$ is topologically isomorphic to $\prod\limits_{\begin{array}{l}_{m\in\mathbb{N}}\\^{p\in\mathbb{P}_G}\end{array}}\langle y_m^{(p)}\rangle$.
\ethm
\begin{proof} Since each group $G_i$ is finite and abelian, it follows that it is a finite sum of finite $p$-groups,
that is $G_i=\bigoplus\limits_{p\in\mathbb{P}_i}(G_i)_p$ and $\mathbb{P}_i=\mathbb{P}_{G_i}$ is finite, $i\in\mathbb{N}$.
Note that $\mathbb{P}_G=\cup\mathbb{P}_i$. Then
$\prod\limits_{i\in\mathbb{N}}G_i\cong\prod\limits_{i\in\mathbb{N}}(\prod\limits_{p\in\mathbb{P}_i}(G_i)_p)\cong\prod
\limits_{p\in\mathbb{P}_G}(\prod\limits_{i\in\mathbb{N}_p}(G_i)_p)$ where $\mathbb{N}_p=\{i\in\mathbb{N}\,:\,G_i\hbox{ has a }p-\hbox{subgroup}\}$.
Consider the embedding
$j:G\hookrightarrow \prod\limits_{p\in\mathbb{P}_G}(\prod\limits_{i\in\mathbb{N}_p}(G_i)_p)$ and the canonical projection
$\pi_p:\prod\limits_{p\in\mathbb{P}_G}(\prod\limits_{i\in\mathbb{N}_p}(G_i)_p)\rightarrow\prod\limits_{i\in\mathbb{N}_p}(G_i)_p$.

Set $G^{(p)}=(\pi_p\circ j)(G)$, that is a compact group. Since $G\cap(\bigoplus\limits_{i\in\mathbb{N}}G_i)$ is dense in $G$,
it follows that
$(\pi_p\circ j)(G\cap(\bigoplus\limits_{i\in\mathbb{N}}G_i))=G^{(p)}\cap(\bigoplus\limits_{i\in \mathbb{N}_p}G_i)_p$ is dense in $G^{(p)}$.
Observe that if $p\in\mathbb{P}_G$ then $G^{(p)}\cap(\bigoplus\limits_{i\in \mathbb{N}_p}G_i)_p=(G\cap(\bigoplus\limits_{i\in \mathbb{N}}G_i))_p$
(otherwise it is the neutral element).
 Applying Theorem \ref{th_products_p}, for each $p\in \mathbb{P}_G$ there is a sequence
 $\{y_m^{(p)}\,:\,m\in\mathbb{N}\}\subseteq (G\cap(\bigoplus\limits_{i\in \mathbb{N}}G_i))_p$
 such that $G^{(p)}\cong\prod\limits_{m\in\mathbb{N}}\langle y_m^{(p)}\rangle$. Then the sequence
 $\{y_m^{(p)}\,:\,m\in\mathbb{N},p\in\mathbb{P}_G\}$ verifies the proof.
\end{proof}
\medskip

The notion of rectangular and weakly rectangular subgroup of an infinite product extend canonically to subgroups of infinite direct sums.
In this direction, we have:

\bthm\label{th_sums}
Let $\mathcal C$  be an order-controllable subgroup of $\bigoplus\limits_{k\in \mathbb{N}} G_k$
such that every group $G_k$ is finite and abelian. Then $\mathcal{C}$ is weakly rectangular. In particular, there is a sequence
$( y_n)\subseteq \mathcal C$ such that
$$\mathcal C\simeq \bigoplus_{n\in \mathbb{N}}\langle y_n \rangle$$
\ethm

\section{Group codes}
According to Forney and Trott \cite{forney_trott:04}, a \emph{group code} is a set of sequences
that has  a group property under a componentwise group operation. In this general setting,
a group code may also be seen as the behavior of a behavioral group system as given by Willens
\cite{willems:86,willems:97}. It is known that many of the fundamental properties of linear codes and systems
depend only on their group structure. In fact, Forney and Trott, loc. cit., obtain purely
algebraic proofs of many of their results. In this section, we follow this approach in order to
apply the results in the preceding sections to obtain further information about the structure of
group codes in very general conditions.

Without loss of generality, assume, from here on, that a \emph{group code} is a subgroup
of a (sequence) group $\mathcal{W}$, called \emph{Laurent group}, that has the generical form
$\mathcal{W}=\mathcal{W}_f\times \mathcal{W}^c$, where $\mathcal{W}_f$  is a direct sum of
abelian groups (locally compact in general) and $\mathcal{W}^c$  a  direct product.
More precisely, let $I\subseteq \mathbb{Z}$ be a countable index set
and let $\{G_k: k\in I\}$ be a set of symbol
groups, a \emph{product sequence space} is a direct product
$$\mathcal{W}^c=\prod _{k\in I} G_k$$
equipped with the canonical product
topology. A \emph{sum sequence space} is a direct sum
$$\mathcal{W}_f=\bigoplus_{k\in I} G_k$$
equipped with the canonical sum (box) topology.
Sequence spaces are often defined to be \emph{Laurent sequences}
$$\mathcal{W}_L=(\bigoplus_{k<0} G_k)\times (\prod_{k\geq 0}G_k).$$
The character group of $\mathcal{W}_L$ is $$\mathcal{W}_{aL}=(\prod_{k<0} \widehat{G_k})\times
(\bigoplus_{k\geq 0}\widehat{G_k}).$$

Thus, a group code  $\mathcal{C}$ is a subgroup of a group
sequence space $\mathcal{W}$ and is equipped with the natural subgroup topology.
Next we recall some basic facts of this theory (cf.
\cite{fagnani:adv97,forney_trott:ie3trans93,forney_trott:04}).
These notions are used in the study of \emph{convolutional codes}
that are well known and used currently in data transmission
(cf. \cite{forney_trott:ie3trans93}).

Let $\mathcal C$ be a group code in the product sequence space $\mathcal W=\prod_{k\in I} G_k$.
According to Fagnani, $\mathcal C$ is called \emph{weakly controllable} if
it is generated by its finite sequences. In other terms, if
$$\overline{\mathcal C}=\overline{\mathcal C\cap \mathcal W_f}.$$

The group code $\mathcal C$ is called \emph{pointwise controllable}
if for all $w_1 , w_2\in  \mathcal C$ and $k\in I$, there exist $L(k)\in \mathbb{N}$ and $w\in \mathcal C$
with $w(i)=w_1(i)\ \forall i<k$, and $w(i)=w_2(i)\ \forall i\geq k+L(k)$.

With the notation introduced above, Let $\mathcal{C}$ be a group code in $\mathcal{W}$. For any $k\in I$
and $L\in \mathbb{N}$, we set
$$\mathcal{C}_k(L):=\{c\in\mathcal{C} :\ \hbox{there exists}\ w\in \mathcal C\ \hbox{with}\ w(i)=0\ \forall i<k\
\hbox{and}\ w(i)=c(i)\ \forall i\geq k+L\}$$

\noindent and $$\mathcal{C}_k:=\bigcup\limits_{L\in\mathbb{N}} \mathcal{C}_k(L).$$

\noindent Obviously $\mathcal{C}_k(1)\subseteq \mathcal{C}_k(2)\subseteq \cdots \mathcal{C}_k(L)\subseteq \cdots\subseteq \mathcal{C}_k$.

We have the following equivalence, whose verification is left to the reader.

\bprp
$\mathcal{C}$ is controllable if and only if $\mathcal{C}=\bigcap\limits_{k\in I}\bigcup\limits_{L\in\mathbb{N}} \mathcal{C}_k(L)$.
\eprp

Given a group code $\mathcal{C}$, the subgroup $$\mathcal{C}_c:=\bigcap\limits_{k\in I}\bigcup\limits_{L\in\mathbb{N}} \mathcal{C}_k(L)$$
is called the \emph{controllable subcode} of $\mathcal{C}$. A code $\mathcal{C}$ is called \emph{uniformly controllable}
when for every $k\in I$, there is $L_k$ such that
$\mathcal{C}=\bigcap\limits_{k\in I} \mathcal{C}_k(L_k)$.
If there is some $L\in\mathbb{N}$ such that $\mathcal{C}_k=\mathcal{C}_k(L)$ for all $k\in I$, it is said that $\mathcal C$ is
\emph{$L$-controllable}. Finally, $\mathcal C$ is \emph{strongly controllable} if it is $L$-controllable for some $L$.
If $\mathcal{C}$ is uniformly controllable and the sequence $(L_k)$ is bounded, then $\mathcal{C}$ is
strongly controllable and the least such $L$ is the \emph{controllability index}
(controller memory) of $\mathcal C$.

Using the same words as in \cite{forney_trott:04}, the core meaning of
``controllable'' is that any code sequence can be reached from any other
code sequence in a finite interval. The following property is clarifying in this regard.
In the sequel $$\mathcal{C}_{:[k,k+L)}=\{w\in\mathcal C : w(j)=0\ \forall j\notin [k,k+L)\}.$$

\bprp
$\mathcal{C}$ is controllable if and only if for any $w\in\mathcal{C}$, there is a sequence
$(L_k)$ contained in $\mathbb{N}$ such that $w\in \sum\limits_{k\in I} \mathcal{C}_{:[k,k+L_k)}$.
\eprp
\begin{proof}
Let $w\in \mathcal{C}$ and let $k_i\in I$ be the first index such that $w(k_1)\not=0$.
Then there is $w_1\in \mathcal{C}$ and $L_1\subset I$ such that $w_1(k_1)=w(k_1)$ and
$w_1(i)=0$ if $i\geq L_1+1$. Take $k_2=L_1+2$ and let $w_2\in \mathcal{C}$ satisfying
$w_2(i)=(w-w_1)(i)$ for all $i<k_2$ and $w_2(i)=0$ for all $i\geq k_2+L_2$.
In general we select $w_n\in \mathcal{C}$ such that $w_n(i)=0$ if $i<k_{n-1}$,
$w_n(i)=(w-w_1-\cdots w_{n-1})(i)$ for all $i<k_3$ and $w_n(i)=0$ if $i\geq k_n+L_n$.
We have that $w=\sum\limits_{n\in \mathbb{N}} w_n$ is the product topology and furthermore
the sum $\sum\limits_{n\in \mathbb{N}} w_n(i)$ is finite for all $i\in I$.
\end{proof}

Analogous notions are defined regarding the observability
of a group code.
The group code $\mathcal C$ is called \emph{weakly observable} if
$$\mathcal C\cap \mathcal W_f= \overline{\mathcal C}\cap \mathcal W_f.$$

Let $\mathcal{C}$ a group code in $\mathcal{W}$, we set
$$(\mathcal{C}_f)_k[L]:=\{c\in\mathcal{W}_f :\ c_{|[k,k+L]}\in \mathcal{C}_{|[k,k+L]}\}.$$

The group code $\mathcal C$ is called \emph{pointwise observable} if
$$\mathcal{C}_f=\bigcap\limits_{k\in I}\bigcap\limits_{L\in\mathbb{N}} (\mathcal{C}_f)_k[L].$$
\medskip

If $\mathcal{C}$ is a group code, then the supergroup
$$\mathcal{C}_{ob}:=\mathcal{C}\cup (\bigcap\limits_{k\in I}\bigcap\limits_{L\in\mathbb{N}} (\mathcal{C}_f)_k[L])$$ is called the \emph{observable supercode}
of $\mathcal{C}$. A code $\mathcal{C}$ is called \emph{uniformly observable}
when for every $k\in I$, there is $L_k$ such that
$\mathcal{C}_f=\bigcap\limits_{k\in I} (\mathcal{C}_f)_k[L_k]$.
If there is some $L\in\mathbb{N}$ such that $\mathcal{C}_f=\bigcap\limits_{k\in I} (\mathcal{C}_f)_k[L]$ for all $k\in I$, it is said that $\mathcal C$ is
\emph{$L$-observable}. Finally, $\mathcal C$ is \emph{strongly observable} if it is $L$-observable for some $L$.
Obviously, if $\mathcal{C}$ is uniformly observable and the sequence $(L_k)$ is bounded, then $\mathcal{C}$ is
strongly observable and the least such $L$ is the \emph{observability index}
(observer memory) of $\mathcal C$.

Recently, Pontryagin duality methods have been applied
systematically  in the study of convolutional abelian group codes.
In this approach, a dual code
$\mathcal C^\perp$ is associated to every group code $\mathcal C$
using Pontryagin-van Kampen duality  in such a way that the
pro\-perties of $\mathcal C$ can be reflected (\emph{dualized}) in
$\mathcal C^\perp$. Along this line,
the following duality theorem provides strong justification for the use
of duality in convolutional group codes (see \cite{forney_trott:04}).
\bthm\label{th:dualcode}\cite{forney_trott:04}
Given dual group codes $\mathcal C$ and $\mathcal C^\perp$,  then
$\mathcal C$ is (resp. weakly, strongly)
controllable if and only if $\mathcal C^\perp$ is
(resp. weakly, strongly) observable, and vice versa.
\ethm


Using duality, we obtain the following additional equivalences (cf. \cite{forney_trott:04}).

\bprp\label{th:dualcode2}
For any group code $\mathcal{C}$ we have
\begin{enumerate}
\item{$(\mathcal{C}_c)^\perp\cong (\mathcal{C}^\perp)_o$.}

\item{$\mathcal{C}$ is controllable if and only if $\mathcal{C}^\perp$ is observable.}

\item{$\mathcal{C}$ is uniformly controllable if and only if $\mathcal{C}^\perp$ is uniformly observable.}
\end{enumerate}
\eprp

Therefore, we can put our attention on the controllability of a group code wlog.
In this direction, the following result was proved in \cite{FHS:2017}.

\bthm[\cite{FHS:2017}]
Let $\mathcal C\leq \prod\limits_{k\in \mathbb{N}} G_k$ be a complete group code
such that every group $G_k$ is finite (discrete). Then the following conditions are equivalent:
\begin{enumerate}
\item{$\mathcal C$ is weakly controllable.}

\item{$\mathcal C$ is controllable.}

\item{$\mathcal C$ is uniformly controllable.}

\end{enumerate}
\ethm
\medskip

In \cite{fagnani:adv97} Fagnani proves that, if $\mathcal{C}$ is a closed, time invariant, group code in $G^\mathbb{Z}$,
with $G$ being a compact group, then the properties of weak controllability,
controllability and strong controllability are equivalent. A different proof of this result follows
easily using the ideas introduced above. Indeed, suppose that $\mathcal{C}$ is a weakly controllable,
compact group code in $\mathcal{W}$. By
Theorem \ref{controllable}, we know that $\mathcal{C}$ is controllable and therefore
$\mathcal{C}=\bigcap\limits_{k\in I}\bigcup\limits_{L\in\mathbb{N}} \mathcal{C}_k(L)$. Suppose, in addition,
that $\mathcal{C}$ is time invariant, then $\mathcal{C}_k(L)=\mathcal{C}_0(L)$ for all $k\in I$. Furthermore,
using Baire category theorem and the compactness of $\mathcal{C}$, it follows that there must
be some $L\in \mathbb{N}$ such that $\mathcal{C}_0=\mathcal{C}_0(L)$, which means that $\mathcal{C}$ is strongly
controllable.

The results formulated above do not hold in general. In fact, an example of a group $H$ that is weakly controllable but not controllable
is provided in \cite{FHS:2017}. Furthermore, using Theorem \ref{th:dualcode}, we obtain that the group $H^\bot$ is is weakly controllable
but not controllable.


As a consequence of the preceding results we obtain the following relation between weakly controllable and
controllable group codes (cf. \cite{FHS:2017}).

\bthm\label{th:code}\label{controllable}
If $\mathcal C$ is a group code in
$$\mathcal{W}=\mathcal{W}_f\times \mathcal{W}^c=(\bigoplus_{i<0} G_i)\times (\prod_{i\geq 0}G_i).$$

\noindent Then the following assertions hold:
\begin{enumerate}
\item[(a)]{If every group $G_i$ is discrete, then $\mathcal C$ is controllable if and only if $\mathcal{C}$ is weakly controllable.}
\item[(b)]{If every group $G_i$ is finite (discrete), then $\mathcal{C}$ is weakly controllable if and only if $\mathcal{C}$ is uniformly controllable.}
\item[(c)]{If every group $G_i$ is a fixed compact group $G$, and $\mathcal{C}$ is a time-invariant, closed subgroup of $\mathcal{W}^c$, then $\mathcal{C}$ is controllable
if and only if $\mathcal{C}$ is strongly controllable.}
\end{enumerate}
\ethm

In case the groups in the family $\{G_i : i\in I\}$ are abelian, Theorem \ref{th:dualcode} yields a similar result for observable group codes,
using Pontryagin duality.
\bthm\label{th:dual_code}
If $\mathcal C$ is a group code in
$$\mathcal{W}=\mathcal{W}_f\times \mathcal{W}^c=(\bigoplus_{i<0} G_i)\times (\prod_{i\geq 0}G_i).$$

\noindent Then the following assertions hold:
\begin{enumerate}
\item[(a)]{If every group $G_i$ is discrete abelian, then $\mathcal C$ is observable if and only if $\mathcal{C}$ is weakly observable.}
\item[(b)]{If every group $G_i$ is finite (discrete) abelian, then $\mathcal{C}$ is weakly observable
if and only if $\mathcal{C}$ is uniformly observable.}
\item[(c)]{If every group $G_i$ is a fixed discrete abelian group $G$, and $\mathcal{C}$ is a time-invariant subgroup of $\mathcal{W}_f$, then $\mathcal{C}$ is observable
if and only if $\mathcal{C}$ is strongly observable.}
\end{enumerate}
\ethm

\section{Conclusion}

To conclude, let us point out that, so far, the applications of Harmonic Analysis and duality methods to the study of group codes have basically reached the abelian case (via Pontryagin duality and Fourier analysis). The non-commutative case has not yet been fully studied, but it can be expected that the application of duality techniques
in the study of non-abelian group codes
could provide some results analogous to those already known in the Abelian case (see the work of Forney and Trott, op.cit). However, the nonabelian duality requires more complicated tools such as Kre\v{i}n algebras, von Neumann algebras, operator spaces, etc.). Therefore, it is first necessary to develop an \emph{appropriate} nonabelian duality that can be applied in a similar way to how it is done in the Abelian case.
\medskip

\textbf{Acknowledgment:} The authors thank Dmitry Shakahmatov for several helpful comments.

\end{document}